\documentclass[12pt,a4paper]{article}
\usepackage{latexsym,amssymb,amsmath,mathrsfs}

\usepackage{sectsty}
\usepackage{color}
\usepackage{bm}
\usepackage[T1]{fontenc}
\usepackage[anchorcolor=blue,%
 bookmarks=true,%
 bookmarksnumbered=true,%
]{hyperref}
\usepackage{cite}
\usepackage{pst-all}
\usepackage{graphicx}

\allsectionsfont{\normalsize\sc\center}

\renewcommand{\theenumi}{{\rm(\roman{enumi})}}

\renewcommand{\theenumii}{{\rm(\alph{enumii})}}

\newtheorem{thm}{Theorem}[section]
\newtheorem{prop}[thm]{Proposition}
\newtheorem{cor}[thm]{Corollary}
\newtheorem{lem}[thm]{Lemma}
\newtheorem{defn}[thm]{Definition}
\newtheorem{remark}[thm]{Remark}
\newtheorem{example}[thm]{Example}

\newenvironment{pf}{\par\begin{trivlist}%
\item[]{\bf Proof.}\ }{\hfill $\square$ \end{trivlist}\par}
\newenvironment{apf}[1]{\par\begin{trivlist}%
\item[]{\bf Proof of #1.}\ }{\hfill $\square$ \end{trivlist}\par}

\makeatletter \@addtoreset{equation}{section} \makeatother

\renewcommand{\P}{\mathbb{P}}

\newcommand{\R}{\mathbb{R}}

\DeclareMathOperator{\Ric}{Ric}

\renewcommand{\d}{\mathrm{d}}
\newcommand{\m}{\mathfrak{m} }

\newcommand{\s}{\mathfrak{s}}

\newcommand{\1}{{\bf 1}}
\newcommand{\eps}{{\varepsilon}}

\title{\large\bf New Laplacian comparison theorem and its applications to diffusion processes  on Riemannian manifolds }
\author{Kazuhiro Kuwae\thanks{Department of Applied Mathematics, Fukuoka University,
Fukuoka 814-0180, Japan ({\sf kuwae@}  {\sf fukuoka-u.ac.jp}) . Supported in part by JSPS Grant-in-Aid for Scientific Research (KAKENHI) 17H02846 and by fund (No.:185001) from the Central Research Institute of Fukuoka University.} 
\ \ and\ \ 
Xiang-Dong Li\thanks{Academy of Mathematics and Systems Science, Chinese Academy of Sciences, 55, Zhongguancun East Road, Beijing, 100190, China, and 
School of Mathematical Science, 
   University of Chinese Academy of Sciences, 
   Beijing 100049, 
   China ({\sf xdli@amt.ac.cn}). 
Supported by National Key R\&D Program of China (No. 2020YFA0712700), NSFC No.~12171458, 11771430, and Key Laboratory RCSDS, CAS, No.~2008DP173182.}
}
\date{}
\begin{document}
\maketitle
\begin{abstract}
Let $L=\Delta-\langle  \nabla\phi, \nabla\cdot\rangle $ be a symmetric diffusion operator with an 
invariant measure $\mu(\d x)=e^{-\phi(x)}\m(\d x)$ on a complete non-compact smooth Riemannian manifold $(M,g)$ with its volume element $\m=\text{\rm vol}_g$, and $\phi\in C^2(M)$ a potential function. 
In this paper, we prove a Laplacian comparison theorem on  weighted complete 
Riemannian manifolds with ${\rm CD}(K, m)$-condition for $m\leq 1$ and a continuous function $K$. 
As consequences, we give the optimal conditions on $m$-Bakry-\'Emery Ricci 
tensor for $m\leq1$ such  that the (weighted) Myers' theorem, Bishop-Gromov volume comparison theorem, 
stochastic completeness and Feller property of $L$-diffusion processes  
hold on weighted complete Riemannian manifolds.  
Some of these results were well-studied for 
$m$-Bakry-\'Emery Ricci curvature for $m\geq n$ (\!\!\cite{Lot,Qian,XDLi05, WeiWylie}) or $m=1$ (\!\!\cite{Wylie:WarpedSplitting, 
WylieYeroshkin}). When $m<1$, our results are new in the literature. 
\end{abstract}


 
\section{Introduction}\label{sec:intro}

 Laplacian comparison theorem is an important result in Riemannian geometry and has many deep applications in geometric analysis on complete Riemannian manifolds, for example,  Myers'  theorem, the Bishop-Gromov volume comparison theorem, the eigenvalue comparison theorem, the Cheeger-Gromoll splitting theorem, the Li-Yau Harnack inequality for positive solution to the heat equation $\partial_t u=\Delta u$ and the upper and lower bounds estimates of heat kernel on complete Riemannian manifolds with natural geometric condition $\Ric\geq K$, where $\Ric$ is the Ricci curvature,  and $K\in \mathbb{R}$.  See S. T. Yau~\cite{Yau75, Yau76, Yau78}, R.  T. Varopoulos \cite{Varo1, Varo2}, Karp and P. Li~\cite{KL},  Li and Yau~\cite{PeterLiYau:Schoedinger}, Schoen and Yau~\cite{SchoenYau:LectDiffGeo} and  the references therein.

On the other hand, by It\^o's theory of stochastic 
 differential equations, the transition probability density of Brownian motion on Riemannian manifolds is the fundamental solution (i.e., the heat kernel) to the heat equation $\partial_t u=\Delta u$.  Due to this important connection between probability theory and geometric analysis,  the Laplacian comparison theorem has also significant  applications in the study of 
   probabilistic aspects of diffusion processes on  complete Riemannian manifolds. In particular,  
     the conservativeness (equivalently, the stochastic completeness) and the Feller property of 
the Brownian motion on manifolds.

 A diffusion process is said to be \emph{conservative} or \emph{stochastically complete} if the associated stochastic process stays in 
the state space forever. This property is equivalent to the 
strong Liouville property for the solution of 
$(L-\lambda)u=0$ for sufficiently large $\lambda>0$, where 
$L$ is the generator associated to the diffusion process.  
More precisely, this means that 
there exists $\lambda_0>0$ such that for all $\lambda\geq\lambda_0$, every non-negative bounded solution of 
$(L-\lambda)u=0$ must be identically zero. 
There are many results on the conservation property
for diffusion processes (see, e.g., 
\cite{Davies:HKconservative, Gaff:Conservative, Grigo:Conservative, HusQin:Conserve, Ichihara:Explosion, Oshima:Conservative, Pang:L1, St:DirI, Takeda:Martingale} and references
therein). In these papers, the conservation property is characterized in terms of the
volume growth of the underlying measure and the growth of the coefficient. The conservativeness of Brownian motions 
on complete Riemannian manifolds has been also 
studied by many authors. 
First, Yau~\cite{Yau76} proved that every complete Riemannian manifold with Ricci curvature bounded from below is stochastically complete. 
Karp and Li~\cite{KL} proved that if the volume of 
the geodesic balls $B_R(x)$ of a complete Riemannian manifold $M$ satisfies $V_R(x)\leq e^{Cr^2}$ for some (and hence all) $x\in M$ and all $r>0$ then $M$ is 
stochastically complete. Li~\cite{PeterLi:Uniqueness} 
proved that if  
${\Ric}(x)\geq- C(1+r_p(x)^2)$ for all $x\in M$ then $M$ is stochastically complete. Li's result can be also considered as a special case of a conservativeness 
criterion due to Varopoulos~\cite{Varo1} and 
Hsu~\cite{Hsu:1989, Hsu:2001}, where they proved that if there exists a non-negative increasing function ${\sf K}(r)$ on $[0,+\infty[$ such that 
\begin{align}
{\Ric}(x)\geq -{\sf K}(r_p(x))\label{eq:Ric}
\end{align}
and 
\begin{align}
\int_{r_0}^{\infty}\frac{\d r}{\sqrt{{\sf K}(r)}}=+\infty\quad \text{ for \ \ some }\quad r_0>0,\label{eq:Inttegrability}
\end{align}
then $M$ is 
stochastically complete.  
So far it is known that the optimal geometric condition for the stochastic completeness of a complete  Riemannian manifold   
is due to Grigor'yan~\cite{Grigo:Conservative} and  in which it was proved that if the volume of geodesic balls of a complete Riemannian manifold $M$ satisfies 
\begin{align}
\int_1^{\infty}\frac{ r\d r}{\log \m(B_r(p))}=+\infty\label{eq:GrigoryanTest}
\end{align} 
for some (and hence all) $p\in M$, then $M$ is stochastically complete. 
Here $\m$ is the volume measure of $(M,g)$. In relation to the Grigoryan's criterion \eqref{eq:GrigoryanTest} for the 
stochastic completeness, Hsu-Qin~\cite{HusQin:Conserve} gave 
a characterization of upper rate function of the process 
in terms of 
 a more relaxed criterion than \eqref{eq:GrigoryanTest}.  
The first example of complete 
but not stochastic complete Riemannian manifold was constructed by Azencott~\cite{Azencott:Behavi}. 
Lyons~\cite{Lyons:Instable} showed that the stochastic completeness is not stable under general quasi-isometric changes of Riemannian metrics.

A diffusion process on $M$  is said to have \emph{Feller property} if its semigroup $P_t$ satisfies $P_t(C_{0}(M))\subset C_{0}(M)$ for each $t>0$ and $\lim_{t\to0}P_tf=f$ for $f\in C_{0}(M)$. 
Here $C_{0}(M)$ denotes the family of continuous functions on $M$ vanishing at infinity. 
A diffusion process on $M$ is said to have \emph{strong Feller property} if its semigroup $P_t$ satisfies $P_t(\mathscr{B}_b(M))\subset C_b(M)$ for each $t>0$. Here 
$\mathscr{B}_b(M)$ (resp.~$C_b(M)$) denotes the family of bounded Borel (resp.~continuous)
functions on $M$. 
By Molchanov~\cite{Molchanov:strongFeller}, 
the semigroup $P_t=e^{tL}$ of the diffusion always 
possesses the strong Feller property.
By Azencott~\cite{Azencott:Behavi}, the semigroup $P_t=e^{tL}$ of the diffusion 
has the Feller property if and only if the following Liouville theorem holds for solutions of $(L-\lambda)u=0$ in 
the exterior region: for any compact set $K\subset M$ and any $\lambda>0$, the minimal positive solution of 
$(L-\lambda)u=0$ on $M\setminus K$ with Dirichlet boundary 
condition $u\equiv 1$ on $\partial K$ must tend to zero at infinity. If ${\bf X}=(\Omega, X_t,\P_x)$ denotes the diffusion on $M$ starting from $x\in M$, then, 
under the strong Feller property of ${\bf X}$, ${\bf X}$ has the Feller property if and only if for each $t>0$ and for all compact set $K\subset M$, 
\begin{align*}
\lim_{d(x,p)\to\infty}\P_x(\sigma_K\leq t)=0,
\end{align*}  
where $p\in M$ and $\sigma_K:=\inf\{t>0\mid X_t\in K\}$ is the first hitting time to $K$ of ${\bf X}$. In \cite{Yau78}, Yau proved that every complete Riemannian manifold with Ricci curvature bounded from below by a negative constant has the  Feller property. We refer the reader to 
Dodziuk~\cite{Dodziuk} for an alternative proof of this result for which one need only to use the maximum principle. By developing Azencott's idea, Hsu~\cite{Hsu:1989,Hsu:2001} proved that if $M$ is complete Riemannian manifold on which there exists a positive increasing 
continuous function ${\sf K}$ on $[0,+\infty[$ satisfying 
\eqref{eq:Ric} and \eqref{eq:Inttegrability}, then the Brownian motion on $M$ has 
the Feller property.

In \cite{Qian}, Z. Qian extended  the standard Laplacian comparison theorem for the usual Laplace-Beltrami operator  to weighted Laplacian (called also Witten Laplacian)  on complete Riemannian manifolds with weighted volume measure and proved an extension of Myer's theorem on 
weighted Riemannian manifolds. In \cite{BQ2}, Bakry and Qian gave a  proof of the  weighted Laplacian comparison theorem without using the Jacobi field theory. In \cite{XDLi05},  the second named author of this paper gave a natural proof of 
Bakry-Qian's  weighted Laplacian comparison theorem  on weighted complete Riemannian manifolds, which is more familiar to readers in geometric analysis, and extended several important 
results in geometric and stochastic analysis to weighted Riemannian manifolds, including  Yau's strong Liouville theorem,  the $L^1$-Liouville theorem, the $L^1$-uniqueness of the solution of the heat equation, the conservativeness and the Feller property for symmetric diffusion processes under the optimal geometry condition on the so-called $m$-dimensional Bakry-\'Emery Ricci curvature, denoted by $\Ric_{m, n}(L)$,  associated with the diffusion operator $L=\Delta-\langle  \nabla\phi,\nabla\cdot\rangle $ 
on an $n$-dimensional complete Riemannian manifold $(M, g)$ with a weighted volume measure $\d\mu=e^{-\phi}\d \m$, where $\phi\in C^2(M)$, $m\geq n$, and
\begin{eqnarray*}
\Ric_{m, n}(L)=\Ric+\nabla^2\phi-{\nabla\phi\otimes \nabla \phi\over m-n}.
\end{eqnarray*}
In particular, the stochastic completeness  and the Feller property for diffusion operators have been proved in \cite{XDLi05}  under the sharp condition on the $m$-dimensional Bakry-Emery Ricci curvature on complete Riemannian manifolds
for $m\geq n$.

 In \cite{XDLi05},  the Li-Yau Harnack inequality has been also proved for positive solutions to the heat equation associated with the weighted Laplacian on weighted complete Riemannian manifolds with the $m$-dimensional Bakry-\'Emery Ricci curvature lower bound condition for $m\geq n$. In \cite{Li12}, the second named author of this paper proved the $W$-entropy formula for the heat equation associated with the Witten Laplacian on weighted complete Riemannian manifolds
 with non-negative $m$-dimensional Bakry-\'Emery Ricci curvature, and pointed out its relationship with the Boltzmann-Shannon entropy and the Li-Yau Hanarck quantity. In \cite{LL15, LL18}, S. Li and X.-D. Li 
 extended the Li-Yau-Hamilton Harnack inequality to positive solutions to the heat equation associated with the Witten Laplacian on weighted  Riemannian manifolds with $\Ric_{m, n}(L)\geq K$, and extended the 
 $W$-entropy formula to the $(K, m)$-super Ricci flows.  They also pointed out the relationship between the Li-Yau-Hamilton Harnack inequality and the $W$-entropy for the heat equation associated with the Witten Laplacian  on weighted 
 complete Riemannian manifolds with the condition $\Ric_{m, n}(L)\geq K$ (equivalently, the CD$(K, m)$-condition), where $m\geq n$ and 
 $K\in \mathbb{R}$ are two constants.  
Indeed,  there have been intensive works  in the literature on the study of geometry and analysis of weighted complete Riemannian manifolds 
with the CD$(K, m)$-condition for $m\geq n$ and $K\in \mathbb{R}$ 
or $K$ being a suitable function on the distance function from a fixed point in manifold. See e.g. \cite{AN, BE1,  BGL_book, BQ1, BQ2, BakryLect1581, BL, FanLiZhang, FLL, XDLi05, Li12, LL15, LL18, Lot, WeiWylie} and reference therein.

 The above mentioned results can be regarded as the natural 
 extensions of the well-known results  in geometric or stochastic analysis for the usual Laplacian on complete Riemannian manifolds with suitable lower bound on the Ricci curvature to the weighted Laplacian on weighted complete Riemannian manifolds with the $m$-dimensional Bakry-\'Emery Ricci curvature $\Ric_{m, n}(L)\geq -K$ for  $m>n$ and for $K\in \mathbb{R}$ or $K$ being a suitable function of the distance function from a fixed point in manifolds. 
 It is natural and interesting to ask the question whether we can extend the well-established results to weighted complete Riemannian manifolds with the geometric condition $\Ric_{m, n}(L)\geq -K$ for $m<n$ and for $K\in \mathbb{R}$ or $K$ being a suitable function of the distance function from a fixed point in manifolds. 
This question was raised by Dominique Bakry to Songzi Li and the second named author of this paper during the Workshop on Stochastic Analysis and Geometry held in AMSS (CAS) at   Beijing in 
November-December 2015. We also noticed that during recent years there are already several papers on the study of geometry on 
weighted Riemannian manifolds with $m$-Bakry-\'Emery Ricci curvature for 
 $m<0$ or $m<1$. For example,  Ohta and Takatsu ~\cite{OhtaTaka} proved  the  $K$-displacement convexity of the R\'enyi type 
 entropy under the  $m$-Bakry-\'Emery Ricci tensor condition ${\Ric}_{m, n}(L)\geq K$,  
 for $m\in]\!-\infty,0\,[\,\cup\,[\,n,+\infty\,[$ and $K\in \mathbb{R}$;  Ohta~\cite{Ohta:KN} and Kolesnikov-Milman~\cite{KolesMilman} simultaneously treated 
the case $m<0$; Ohta~\cite{Ohta:KN} extended the Bochner inequality, eigenvalue estimates, and the Brunn-Minkowski inequality 
under the lower bound for 
${\Ric}_{m,n}(L)$ with $m<0$; 
Kolesnikov-Milman~\cite{KolesMilman} also proved the Poincar\'e and the Brunn-Minkowski inequalities for manifolds with boundary 
under the lower bound for 
${\Ric}_{m,n}(L)$ with $m<0$; Ohta ~\cite[Theorem~4.10]{Ohta:KN} also proved that the lower bound of ${\Ric}_{m,n}(L)(x)$ with 
 $m<0$ is equivalent to the curvature dimension condition in terms of 
 mass transport theory as defined by Lott-Villani~\cite{LV2} and Sturm~\cite{St:geomI, St:geomII}.  In 
 \cite{Wylie:WarpedSplitting}, Wylie  proved a generalization of the Cheeger-Gromoll splitting theorem 
 under the ${\rm CD}(0, 1)$-condition. In \cite{WylieYeroshkin}, 
W.~Wylie and D.~Yeroshkin proved  a Bishop-Gromov volume comparison theorem, a Laplacian comparison theorem, Myers' theorem and Cheng's maximal diameter theorem on manifolds with 
$m$-Bakry-\'Emery Ricci curvature condition for $m=1$.  Recently,  Milman~\cite{Milman17} extended the Heintze-Karcher Theorem, isoperimetric inequality, and functional inequalities under the lower bound for ${\Ric}_{m,n}(L)(x)$ with $m < 1$.  

It is important to know whether one can establish the Laplacian comparison theorem and to develop geometric and stochastic analysis 
on weighted complete Riemannian manifolds with 
$\Ric_{m, n}(L)\geq K$  for $m\leq 1$ and $K\in \mathbb{R}$ or $K$ being a suitable function on $M$, because the condition $\Ric_{m, n}(L)\geq K$  for $m\leq 1$ weaker than the condition $\Ric_{m, n}(L)\geq K$  for $m\in[n,+\infty]$ and the results on stochastic or geometric analysis under $\Ric_{m, n}(L)\geq K$  for $m\leq 1$ are few at present. This is one of the reason why the establishment of the Laplacian comparison theorem is so important.
In this paper, we prove such a new Laplacian comparison theorem for 
 $m<1$ and for $K$ being a function depending on a distance 
 function on $M$. As consequences, we give the optimal conditions on the $m$-Bakry-\'Emery Ricci 
tensor for $m\leq1$ so that  the Bishop-Gromov volume comparison theorem, 
and Myers' theorem 
hold on weighted complete Riemannian manifolds.  Moreover, we use the new Laplacian comparison theorem  to establish the stochastic completeness and the Feller property for the $L$-diffusion processes on complete Riemannian manifolds with the optimal geometric condition on $\Ric_{m, n}(L)\geq K$  for $m\leq 1$ and $K\in \mathbb{R}$ or $K$ being a suitable function on $M$. 
 As far as we know, when $m<1$,  our results are new in the literature. 
 
Finally, to end this Introduction,  we would like to point out that, once the new  Laplacian comparison theorem is established, it will be 
possible to develop further study of geometric analysis on weighted complete Riemannian manifolds with the ${\rm CD}(K, m)$-condition for $m\leq 1$ and $K\in \mathbb{R}$. 
In a  joint project with Songzi Li, we will study the gradient estimates and Liouville theorems  for $L$-harmonic functions and the 
Li-Yau Harnack inequality for positive solutions to the heat equation $\partial_t u=Lu$ on weighted complete Riemannian manifolds the ${\rm CD}(K, m)$-condition for 
 $m\leq1$ or $m\leq 0$ and $K\in \mathbb{R}$. 

\section{Main result}

Let $(M,g)$ be an $n$-dimensional smooth complete Riemannian manifold, and 
$\phi\in C^2(M)$  a potential function. Throughout this paper, we assume that the manifold $M$ has no boundary. 
We consider a diffusion operator $L:=\Delta- \langle  \nabla\phi, \nabla\cdot\rangle $, which is symmetric with respect to the invariant measure $\mu(\d x)=e^{-\phi(x)}\m(\d x)$, where $\m:=\text{\rm vol}_g$ is the volume element of $(M,g)$.  In \cite{Li12, LL18},  $L$ is called the Witten Laplacian on $(M, g, \phi)$. 

For any constant $m\in]-\infty,+\infty]$, we introduce the symmetric $2$-tensor 
\begin{align*}
{\Ric}_{m,n}(L)(x)={\Ric}(x)+\nabla^2\phi(x)-\frac{\nabla\phi(x)\otimes\nabla\phi(x)}{m-n},\quad x\in M,
\end{align*} 
and call it the \emph{$m$-Bakry-\'Emery Ricci curvature} of the diffusion operator 
$L$.  For any $m\in]-\infty,+\infty]$ and a continuos function $K:M\to\mathbb{R}$, we call $(M, g, \phi)$ or $L$ satisfies the ${\rm CD}(K, m)$-condition  if 
\begin{align*}
{\Ric}_{m, n}(L)(x)\geq K(x)\quad \text{ for all }\quad x\in M.
\end{align*} 
When $m=n$, we always assume that $\phi$ is a constant so that 
${\Ric}_{n,n}(L)={\Ric}$. 
When $m\geq n$, $m$ is regarded as an upper bound of the dimension of the diffusion operator $L$. 
Throughout this paper, we focus on the case $m\leq1$ and assume $n>1$ if $m=1$ and $\phi$ is not a constant (i.e., $\phi$ is a constant and $L=\Delta$ if $m=n=1$). 
  
\bigskip 

For two points 
$p,q\in M$, the \lq\lq re-parametrized distance\rq\rq\, between $p$ and $q$ is defined  to be 
\begin{align}
s(p,q):=\inf\left\{\left.
\int_0^{r_p(q)}e^{-2\frac{\phi(\gamma_t)}{n-m}}\d t\;\;\right|\left.\begin{array}{ll}
&\!\!\!\!\gamma:\text{unit speed geodesic}\\
&\!\!\!\!{\gamma}(0)=p, {\gamma}(r_p(q))=q
\end{array}\right.\right\}.\label{eq:reparametrizeddistance}
\end{align} 
If $(M,g)$ is complete, then $s(p,q)$ is finite  and 
well-defined from the basic properties of Riemannian geodesics. 
Let $s_p(\cdot):=s(p,\cdot)$. If $q$ is not a cut point of $p$, then there is  a unique minimal geodesic from $p$ to $q$ and $s_p$ is smooth in a neighborhood of $q$ as 
can be computed by pulling the function back by the exponential map 
at $p$. Note that $s(p,q)\geq0$, it is zero if and only if $p=q$, 
and $s(p,q)=s(q,p)$ holds. However, $s(p,q)$ does not define a distance since it does not satisfy the triangle inequality.

\begin{defn}\label{df:phicompletness}
{\rm Let $(M,g)$ be an $n$-dimensional complete Riemannian manifold and $\phi\in C^2(M)$. Take $p\in M$.  Then we say that $(M,g,\phi)$ is \emph{$(\phi,m)$-complete at $p$} if 
\begin{align}
\varlimsup_{r\to+\infty}\inf
\left\{\left. \int_0^re^{-2\frac{\phi(\gamma_t)}{n-m}}\d t\;\right| 
\left.\begin{array}{ll}
&\!\!\!\!\gamma:\text{unit speed geodesic}\\
&\!\!\!\!{\gamma}(0)=p, L(\gamma)=r
\end{array}\right.
\right\}=+\infty.\label{eq:phimcomplete}
\end{align}
We say that  $(M,g,\phi)$ is \emph{$(\phi,m)$-complete} if it is $(\phi,m)$-complete at all $p\in M$.}  
\end{defn}

\begin{remark}
{\rm 
\begin{enumerate}
\item If $(M,g,\phi)$ is $(\phi,m)$-complete at $p$, then 
$\int_0^{\infty}e^{-2\frac{\phi(\gamma_t)}{n-m}}\d t=+\infty$ always holds for any unit speed geodesic $\gamma$ 
with $\gamma_0=p$. 
\item If $\phi$ is upper bounded,  then $(M,g,\phi)$ is $(\phi,m)$-complete. In particular, if $M$ is compact, then 
 $(M,g,\phi)$ is $(\phi,m)$-complete. 
 \item If $C\in]0,(n-m)/2]$ and $\sup_{B_r(p)}\phi\leq C\log(1+r)$ holds 
 for all $r\in[r_0,+\infty[$ with
 some $r_0>0$, then $(M,g,\phi)$ is $(\phi,m)$-complete at $p$. In fact,  
 \begin{align*}
 \inf_{L(\gamma)=r}\int_0^re^{-2\frac{\phi(\gamma_t)}{n-m}}\d t\geq  \int_{r_0}^re^{-\frac{2C\log (1+t)}{n-m}}\d t=\int_{r_0}^{r}\frac{\d t}{(1+t)^{2C/(n-m)}}\to +\infty \quad \text{ as }\quad r\to\infty.
 \end{align*}
\item The $(\phi,1)$-completeness defined in \cite[Definition~6.2]{Wylie:WarpedSplitting} implies the $(\phi,m)$-completeness  
 provided $\phi$ is non-negative. 
The converse also holds under the non-positivity of $\phi$.  
\end{enumerate}
}
\end{remark}
\begin{lem}\label{lem:phimcomplete}
Let $(M,g)$ be an $n$-dimensional complete Riemannian manifold and $\phi\in C^2(M)$. Take $p\in M$ and suppose that \eqref{eq:phimcomplete} holds. Then, for any sequence $\{q_i\}$ in $M$  such that  $d(p,q_i)\to+\infty$ as $i\to+\infty$, 
$s(p,q_{i})\to+\infty$ as $i\to+\infty$.
\end{lem}
\begin{pf}
The proof is similar to that of \cite[Proposition~3.4]{WylieYeroshkin}. We omit it.  
\end{pf}

\subsection{Laplacian Comparison} 
Let $\kappa:[0,+\infty[\to\R$ be a continuous function and ${\sf a}_{\kappa}$ the unique solution defined on the maximal interval $]0,\delta_{\kappa}[$ for $\delta_{\kappa}\in]0,+\infty]$ of the following Riccati equation 
\begin{align}
-\frac{\d {\sf a}_{\kappa}}{\d s}(s)=\kappa(s)+{\sf a}_{\kappa}(s)^2\label{eq:RiccatiEq}
\end{align}
with the boundary conditions 
\begin{align}
\lim_{s\downarrow 0}s\, {\sf a}_{\kappa}(s)=1,\label{eq:BdryCond}
\end{align}
and 
\begin{align}
\lim_{s\uparrow \delta_{\kappa}}(s-\delta_{\kappa})\, {\sf a}_{\kappa}(s)=1\label{eq:BdryCondStrict*}
\end{align}
under $\delta_{\kappa}<\infty$. 
\eqref{eq:BdryCond} yields 
$\lim_{s\downarrow0}{\sf a}_{\kappa}(s)=+\infty$.
If $\delta_{\kappa}<\infty$, from \eqref{eq:BdryCondStrict*}, $\delta_{\kappa}$ is the explosion time of ${\sf a}_{\kappa}$ in the sense of  
$\lim_{s\uparrow\delta_{\kappa}}{\sf a}_{\kappa}(s)=-\infty$. 
In fact, ${\sf a}_{\kappa}$ is given by 
${\sf a}_{\kappa}(s)={\s_{\kappa}'(s)}/{\s_{\kappa}(s)}$, where $\s_{\kappa}$ is a unique solution of 
Jacobi equation $\s_{\kappa}''(s)+\kappa(s)\s_{\kappa}(s)=0$ with $\s_{\kappa}(0)=0$, $\s_{\kappa}'(0)=1$, and 
$\delta_{\kappa}=\inf\{s>0\mid \s_{\kappa}(s)=0\}$. We write ${\sf a}_{\kappa}(s)=\cot_{\kappa}(s)$. 
We can deduce 
$\s_{\kappa}(\delta_{\kappa})=0$ and $\s_{\kappa}'(\delta_{\kappa})<0$, hence $\cot_{\kappa}(\delta_{\kappa})=-\infty$ provided $\delta_{\kappa}<\infty$. 
If $\kappa$ is a real constant, then 
$\delta_{\kappa}=\pi/\sqrt{\kappa^+}\leq+\infty$. 
Take $m\in ]-\infty,1\,]$ and  
set $m_{\kappa}(s):=(n-m)\cot_{\kappa}(s)$. 
Then \eqref{eq:RiccatiEq} is equivalent to 
\begin{align}
-\frac{\d m_{\kappa}}{\d s}(s)=(n-m)\kappa(s)+\frac{m_{\kappa}(s)^2}{n-m},\label{eq:RiccatiEqM}
\end{align}
and \eqref{eq:BdryCond} (resp.~\eqref{eq:BdryCondStrict*}) is equivalent to 
$\lim_{s\downarrow 0}s\, m_{\kappa}(s)=n-m$ (resp.~$\lim_{s\uparrow \delta_{\kappa}}(s-\delta_{\kappa})\, m_{\kappa}(s)=n-m$ under $\delta_{\kappa}<\infty$).

Our first result is the following Laplacian comparison theorem on 
weighted  complete Riemannian manifolds 
having lower bound $K(x)=(n-m)\kappa(s_p(x))e^{-\frac{4\phi(x)}{n-m}}$  
for $m$-Bakry-\'Emery Ricci tensor ${\Ric}_{m,n}(L)_x$ on the direction $\nabla r_p$ with $m\leq1$. 

\begin{thm}[Laplacian Comparison Theorem]\label{thm:GlobalLapComp}
Suppose that $(M,g)$ is an $n$-dimen\-sional complete smooth Riemannian manifold and $\phi\in C^2(M)$. 
Take $p\in M$ and $R\in]0,+\infty]$.  
Suppose that 
\begin{align}
{\Ric}_{m,n}(L)_x(\nabla r_p,\nabla r_p)\geq(n-m)\kappa(s_p(x)) e^{-\frac{4\phi(x)}{n-m}}\label{eq:RiciLowBdd}
\end{align} 
holds under $r_p(x)<R$ with $x\in ({\rm Cut}(p)\cup\{p\})^c$.  
Then  
\begin{align}
(L r_p)(x)\leq e^{-\frac{2\phi(x)}{n-m}}m_{\kappa}(s_p(x)).\label{eq:GloLapComp}
\end{align}
\end{thm}

\subsection{Geometric consequences}

\begin{thm}[Weighted Myers' Theorem]\label{thm:WeightedMyers}
Let $(M,g)$  be an $n$-dimensional complete Riemannian manifold and a function $\phi\in C^2(M)$. 
Take $p\in M$.  
Assume that \eqref{eq:RiciLowBdd} 
holds for all $x\in ({\rm Cut}(p)\cup\{p\})^c$ and $\delta_{\kappa}<\infty$.  
  Then $s(p,q)\leq \delta_{\kappa}$ for all $q\in M$. If further assume the 
  $(\phi,m)$-completeness at $p$, then $M$ is compact. 
\end{thm}

\begin{remark}
{\rm 
\begin{enumerate}
\item Theorem~\ref{thm:WeightedMyers} 
generalizes \cite[Theorem~2.2 and Corollary~2.3]{WylieYeroshkin}.
\item Since upper boundedness of $\phi$ implies the $(\phi,m)$-completeness, 
Theorem~\ref{thm:WeightedMyers} 
implies the compactness of $M$ 
 if $\delta_{\kappa}<\infty$, \eqref{eq:RiciLowBdd} holds for all $x\notin({\rm Cut}(p)\cup\{p\})^c$, and $\phi$ is upper bounded. 
\end{enumerate}
}
\end{remark}

Next  we will give 
the Bishop-Gromov type volume comparison. 
This is described for $\mu(A)=\int_Ae^{-\phi(x)}\m(\d x)$ of metric annuli $A(p,r_0,r_1):=\{x\in M\mid r_0\leq r_p(x)\leq r_1\}$. The comparison in this case will be in terms of the quantities
\begin{align*}
\hspace{-0.5cm}\overline{\nu}_p(\kappa,r_0,r_1):=\int_{r_0}^{r_1}\int_{\mathbb{S}^{n-1}}\!\!\s_{\kappa}^{n-m}\left({\color{black}{\sup_{\eta
}s_p(r,\eta)}}
\right)\d r\d\theta,
\quad 
\underline{\nu}_p(\kappa,r_0,r_1):=\int_{r_0}^{r_1}\int_{\mathbb{S}^{n-1}}\!\!\s_{\kappa}^{n-m}\left(
{\color{black}{\inf_{\eta
}s_p(r,\eta)}}
\right)\d r\d\theta.
\end{align*} 
\begin{thm}[Bishop-Gromov Volume Comparison]\label{thm:BGVol}
Take $p\in M$ {\color{black}{and $R\in]0,+\infty]$}}. 
Suppose that $(M,g)$ is an $n$-dimensional complete smooth Riemannian manifold and a function $\phi\in C^2(M)$. 
Let $\kappa:[0,+\infty[\to\R$ be a continuous function.  Assume that \eqref{eq:RiciLowBdd} 
holds {\color{black}{for $r_p(x)<R$ with }}$x\in ({\rm Cut}(p)\cup\{p\}	)^c$. 
Suppose that $0\leq r_0< r_a\leq r_1$ and $0\leq r_0\leq r_b<r_1$. Then
\begin{align}
\frac{\mu(A(p,r_b,r_1))}{\mu(A(p,r_0,r_a))}\leq \frac{\overline{\nu}_p(\kappa,r_b,r_1)}{\underline{\nu}_p(\kappa,r_0,r_a)}\label{eq:BGAnnuliUpLow}
\end{align}
holds {\color{black}{for $r_1<R$}}. 
\end{thm}

\subsection{Probabilistic consequences}
For a fixed $p\in M$, we set  
$\underline{\phi}_p(r):=\inf_{B_r(p)}\phi$ and 
$\overline{\phi}_p(r):=\sup_{B_r(p)}\phi$
for $r>0$. Then $\underline{\phi}_p(r)\leq\phi(p)\leq 
\overline{\phi}_p(r)$ and $\lim_{r\to0}\overline{\phi}_p(r)=\lim_{r\to0}\underline{\phi}_p(r)=\phi(p)$.
It is easy to see 
that for any $y\in M$ with $d(p,y)\leq d(x,p)$
\begin{align*}
\underline{\phi}_p(r_p(x))\leq\phi(y)\leq\overline{\phi}_p(r_p(x)),
\end{align*}
and for each $s>0$
\begin{align}
\underline{\phi}_p(r+s)\leq\underline{\phi}_q(r)
\leq \overline{\phi}_q(r)\leq \overline{\phi}_p(r+s)
\quad\text{ for \ \ any }\quad q\in M\quad\text{ with }\quad d(p,q)\leq s.\label{eq:phipIneQ}
\end{align} 
By \eqref{eq:reparametrizeddistance}, 
we see 
\begin{align}
e^{-\frac{2\overline{\phi}_p(r_p(x))}{n-m}}r_p(x)\leq 
s_p(x)\leq e^{-\frac{2\underline{\phi}_p(r_p(x))}{n-m}}r_p(x).\label{eq:spIneQrp}
\end{align}

Let ${\sf K}(s)$ be a non-negative continuous non-decreasing function on $[0,+\infty[$. 
We consider the following two conditions: 
\begin{align}
\hspace{-2.5cm}{\bf K}(p):\qquad\qquad\int_{r_0}^{\infty}\frac{\d r}{\sqrt{{\sf K}\left(e^{-\frac{2\underline{\phi}_p(r)}{n-m}}r\right)}e^{-\frac{2\underline{\phi}_p(r)}{n-m}}}=+\infty
\quad \text{ for \ \ some }\quad r_0>0.
\label{eq:integrability} \\
\hspace{-1cm}\overline{\bf K}(p):\quad\quad
\int_{r_0}^{\infty}\frac{\d r}{\sqrt{{\sf K}\left(e^{-\frac{2\underline{\phi}_p(r)}{n-m}}r\right)}e^{\frac{(2\overline{\phi}_p(r)-2\underline{\phi}_p(r))}{n-m}}}=+\infty
\quad \text{ for \ \ some }\quad r_0>0.
\label{eq:upintegrability}
\end{align}

\begin{remark}\label{rem:Integrability}
{\rm  
Since $\overline{\phi}_p(r)\geq\phi(p)\geq \underline{\phi}_p(r)$, \eqref{eq:integrability} (resp.~\eqref{eq:upintegrability})  always implies \eqref{eq:Inttegrability} (resp.~\eqref{eq:integrability}), and the converse holds 
provided $\phi$ is lower (resp.~upper) bounded.
}
\end{remark}

\begin{prop}\label{pr:Independence}
The conditions ${\bf K}(p)$ and $\overline{\bf K}(p)$
 are independent of the choice of $p\in M$ respectively. 
\end{prop}
\begin{pf}
Suppose that ${\bf K}(p)$ holds. 
Take another point $q\in M$ and $s>0$ with $d(p,q)\leq s$. 
Then $\underline{\phi}_p(r+s)\leq\underline{\phi}_q(r)$. We may assume $r_0>s$. 
Then 
\begin{align*}
\int_{r_0-s}^{\infty}\frac{\d r}{\sqrt{{\sf K}\left(e^{-\frac{2\underline{\phi}_q(r)}{n-m}}r\right)}e^{-\frac{2\underline{\phi}_q(r)}{n-m}}}
&
\geq
\int_{r_0}^{\infty}\frac{\d r}{\sqrt{{\sf K}\left(e^{-\frac{2\underline{\phi}_p(r)}{n-m}}r\right)}e^{-\frac{2\underline{\phi}_p(r)}{n-m}}}
=+\infty.
\end{align*}
Therefore ${\bf K}(q)$ holds. The proof of the independence of $p\in  M$ for $\overline{\bf K}(p)$ is similar.
\end{pf}

\begin{thm}[Conservativeness of $L$-diffusion]\label{thm:conservative}
Take $p\in M$. 
Let $(M,g)$ be an $n$-dimensional complete Riemannian manifold without boundary and a function $\phi\in C^2(M)$. 
Let ${\sf K}(r)$ be a non-negative continuous non-decreasing  function on $[0,+\infty[$.
Assume that 
\begin{align}
\text{\rm Ric}_{m,n}(L)(x)\geq -{\sf K}(s_p(x))e^{-\frac{4\phi(x)}{n-m}}\quad \text{ for \ \ any }\quad x\in M\label{eq:RicciBound*}
\end{align}
holds. Suppose one of the following: 
\begin{enumerate}
\item ${\sf K}\not\equiv0$ satisfies ${\bf K}(p)$. 
\item  ${\sf K}\equiv0$
and 
\begin{align}
\int_{r_0}^{\infty}
\exp\left({\frac{2\underline{\phi}_p(r)}{n-m}
}\right)\d r
=+\infty\quad \text{ for \ \ some }\quad r_0>0.\label{eq:asymptoticphi}
\end{align}
\end{enumerate}
\noindent Then the heat semigroup  $P_t=e^{tL}$ is conservative, i.e., 
\begin{align*}
P_t1(x)=1, \quad t>0, \quad x\in M.
\end{align*}
Equivalently, there exists some $\lambda_0>0$ such that for all $\lambda\geq\lambda_0$, every non-negative bounded solution of $(L-\lambda)u=0$ must be identically zero. 
Moreover, under (ii), if we assume $n\leq m+1$ and the lower boundedness of $\phi$, then the heat semigroup  $P_t=e^{tL}$ is recurrent, i.e., 
any bounded $L$-subharmonic function is constant. 
\end{thm}

\begin{thm}[Feller Property of $L$-diffusion]\label{thm:Fellerproperty}
Take $p\in M$. 
Let $(M,g)$ be an $n$-dimensional complete Riemannian manifold without boundary and a function $\phi\in C^2(M)$. 
Let ${\sf K}(s)$ be a non-negative continuous non-decreasing  function on $[0,+\infty[$. 
Assume that 
\begin{align}
\inf\left\{\left. \text{\rm Ric}_{m,n}(L)(x)e^{\frac{4\phi(x)}{n-m}}\;\right|\; s_p(x)=s\; \right\}\geq -{\sf K}(s)
\label{eq:RicciBound}
\end{align}
holds for any $s\in [0,\infty[$.
Suppose one of the following: 
\begin{enumerate}
\item ${\sf K}\not\equiv0$ satisfies $\overline{\bf K}(p)$.
\item ${\sf K}\equiv0$ and 
\begin{align}
\int_{r_0}^{\infty}
\exp\left({-\frac{2\overline{\phi}_p(r)-2\underline{\phi}_p(r)}{n-m}
}\right)\d r
=+\infty\quad \text{ for \ \ some }\quad r_0>0.\label{eq:integarbilityK0}
\end{align}
\end{enumerate}
\noindent Then the heat semigroup $P_t=e^{tL}$ has the Feller property. 
Equivalently, for any $\lambda>0$ and any compact subset $K\subset M$, the minimal positive solution 
of $(L-\lambda)u=0$ on $M\setminus K$ with the Dirichlet boundary condition $u\equiv 1$ on $\partial K$ must be zero at infinity. 
\end{thm}
\begin{remark}\label{rem:Integrability*}
{\rm 
\begin{enumerate}
\item Note that \eqref{eq:RicciBound} implies \eqref{eq:RicciBound*}.  
\item Since $\overline{\phi}_p(r)\geq \phi(p)$, \eqref{eq:integarbilityK0} implies \eqref{eq:asymptoticphi}. 
The converse holds if $\phi$ is upper bounded. 
If 
\begin{align}
\varliminf_{r\to\infty}\frac{\;\underline{\phi}_p(r)\;}{\log r}>-\infty,\label{eq:asymptoticphi*}
\end{align}
then \eqref{eq:asymptoticphi} is satisfied. 
\item As in Proposition~\ref{pr:Independence}, 
 \eqref{eq:asymptoticphi} and \eqref{eq:integarbilityK0} are 
 independent of the choice of $p\in M$ respectively.  
When ${\sf K}\equiv K>0$, \eqref{eq:integrability} (resp.~\eqref{eq:upintegrability})  
is equivalent to \eqref{eq:asymptoticphi} (resp.~\eqref{eq:integarbilityK0}).\end{enumerate}
}
\end{remark}

\begin{example}
{\rm Take $\eps\in[0,1]$ and  $\delta\in[\frac{2\eps}{1+\eps},1]$. 
Let $K$ be a positive constant and set 
${\sf K}(r)=Kr^{2(1-\delta)}$. Then \eqref{eq:Inttegrability} always holds. 
We set a non-positive $C^2$-function  
$\phi(x)=-\frac{(n-m)\eps}{4} \log(1+ r_p^2(x))$.  
Then $0\geq\overline{\phi}_p(r)\geq\underline{\phi}_p(r)\geq-\frac{(n-m)\eps}{4} \log(1+ r^2)$, hence 
\eqref{eq:asymptoticphi} and \eqref{eq:integarbilityK0}
hold. Moreover,  
\eqref{eq:integrability} and \eqref{eq:upintegrability}
 also hold in view of 
Remark~\ref{rem:Integrability*}(ii). 
In particular, the case for $\eps=\delta=1$ tells us the following: 
If ${\rm Ric}_{m,n}(L)(x)\geq -Ke^{-\frac{4\phi(x)}{n-m}}=-K(1+r_p(x)^2)$ holds for all $x\in M$ under $m\in]-\infty,1]$, then \eqref{eq:RicciBound} holds, consequently 
the heat semigroup $P_t=e^{tL}$ is conservative and has the Feller property.
}
\end{example}

\begin{example}
{\rm 
We consider the case $(M,g)=(\R^n, g_{\tiny{\rm Euc}})$. 
In this case, ${\rm Ric}_{m,n}(L)=(n-m)e^{-\frac{\phi}{n-m}}
{\rm Hess}\left(e^{\frac{\phi}{n-m}}\right)$. 
Let $K$ be a non-negative constant. 
We consider the following condition:
\begin{enumerate}
\item[{\bf (A)}] \qquad ${\rm Hess}\left(e^{\frac{\phi}{n-m}}\right)+\frac{K}{n-m}e^{-\frac{3\phi}{n-m}}E_n\geq O$. 
\end{enumerate}
Here $E_n$ (resp.~$O$) denotes the $(n,n)$-identity 
(resp.~$(n,n)$-zero)  
matrix. Under ${\bf (A)}$, we have ${\rm Ric}_{m,n}(L)\geq -Ke^{-\frac{4\phi}{n-m}}$ on $\R^n$, that is, \eqref{eq:RicciBound} is satisfied. 
In this case, \eqref{eq:asymptoticphi}, 
\eqref{eq:integarbilityK0} and 
\eqref{eq:asymptoticphi*} 
are satisfied for 
$\phi(x)= -\frac{n-m}{4}\log(1+ |x-p|^2)$ with some $p\in \R^n$. 
So \eqref{eq:integrability} and \eqref{eq:upintegrability} 
are satisfied for such $\phi$ and ${\sf K}\equiv K>0$. Let us consider the case $K=0$. Then   
the heat semigroup $P_t=e^{tL}$ is conservative and has the Feller property under the convexity of $e^{\frac{\phi}{n-m}}$ with 
$m\leq1$. In particular, the convexity of $\phi$ also yields  
the same conclusion, because ${\rm Hess}(\phi)\geq0$ implies ${\rm Hess}(\phi)\geq \frac{\nabla\phi\otimes\nabla\phi}{m-n}$ under $m\leq1$, which means 
${\rm Hess}\left(e^{\frac{\phi}{n-m}}\right)\geq0$.  
Now we consider the case $K>0$. The assumption ${\bf (A)}$ is satisfied for 
 $\phi(x):=-\frac{n-m}{4}\log(1+ |x|^2)$ provided 
 $K\geq\frac{n-m}{2}$. 
Indeed, 
$f(x):=e^{\frac{\phi(x)}{n-m}}=(1+|x|^2)^{-\frac{1}{4}}$ implies  
\begin{align*}
{\rm Hess}\,f(x)&+\frac{K}{n-m}f^{-3}(x)E_n\\
&\geq {\rm Hess}\,f(x)+\frac{K}{n-m}f(x)E_n\\
&=(1+|x|^2)^{-\frac94}
\left(\frac54 x_ix_j-\frac12(1+|x|^2)\delta_{ij}
+\frac{K}{n-m}(1+|x|^2)^2\delta_{ij}\right)_{ij}\\
&\geq \frac54(1+|x|^2)^{-\frac94}(x_ix_j)_{ij}
\geq O
\end{align*}
under $K\geq\frac{n-m}{2}$. 
Therefore 
\eqref{eq:RicciBound} holds under $K\geq\frac{n-m}{2}$. 
Moreover, $\phi$ 
satisfies ${\sf K}(0)$. Hence 
the heat semigroup $P_t=e^{tL}$ is conservative and has the Feller property under $K\geq\frac{n-m}{2}$ 
with $m<1$ for $\phi(x)=-\frac{n-m}{4}\log(1+|x|^2)$. 
}
\end{example}
\section{Proof of Theorem~\ref{thm:GlobalLapComp}}

\subsection{Volume element comparison}
 
Let $p\in M$ and let $(r,\theta)$, $r>0$, $\theta\in \mathbb{S}^{n-1}$ be exponential 
polar coordinates (for the metric $g$) around $p$ which are defined on 
a maximal  star shaped domain in $T_pM$ called the 
\emph{segment domain}. Write the volume element $\d\m=J(r,\theta)\d r\land\d\theta$. 

Let $s_p(\cdot)$ be the re-parametrized distance function defined above. 
Inside the segment domain, 
$s_p$ has the simple formula
\begin{align*}
s_p(r,\theta)=\int_0^re^{-\frac{2\phi(t,\theta)}{n-m}}\d t.
\end{align*}
Therefore, $s$ is a smooth function in the segment domain with the property that 
$\frac{\partial s}{\partial r}=e^{-\frac{2\phi}{n-m}}$. We can then also take $(r,\theta)$ 
to be coordinates which are also valid for the entire segment theorem. We can not control the derivative of $s$ in directions tangent to the sphere, so the new 
$(s,\theta)$ coordinates are \emph{not} orthogonal as in the case for geodesic polar coordinates. However, this is not the issue when we compute volumes as 
\begin{align}
e^{-\frac{2\phi}{n-m}}\d\mu=e^{-\frac{n-m+2}{n-m}\phi}J(r,\theta)\d r\land\d\theta=e^{-\phi}J(r,\theta)\d s\land\d\theta.\label{eq:conformal}
\end{align} 

In geodesic polar coordinates $\frac{\d}{\d s}$ has the expression $\frac{\d }{\d s}=e^{\frac{2\phi}{n-m}}\frac{\partial}{\partial r}$. 
Note that it is not the same as $\frac{\partial}{\partial s}$ in $(s,\theta)$ coordinates. 

Recall that for a Riemannian manifold $\frac{\d}{\d r}\log J(r,\theta)=\Delta r_p$, where 
$\Delta r_p$ is the standard Laplacian acting on the distance function $r_p$ from the point $p$.  \eqref{eq:conformal} indicates we should consider the quantity 
\begin{align}
\frac{\d}{\d s}\log(e^{-\phi}J(r,\theta))=e^{\frac{2\phi}{n-m}}\left(\Delta r_p-\langle  \nabla \phi,\nabla r_p\rangle  \right).\label{eq:quantity}
\end{align}
We thus recover the Witten Laplacian $Lu:=\Delta u-\langle  \nabla \phi, \nabla u\rangle $.  Letting $\lambda=e^{\frac{2\phi}{n-m}}Lr_p$, we find that 
$\lambda$ satisfies the Riccati differential inequality in terms of the parameter $s$.

\begin{lem}\label{lem:RiccatiDiferIneq}
Let $\gamma$ be a unit speed minimal geodesic with $\gamma_0=p$ and $\dot{\gamma}_0=\theta$. 
Let $s$ be the 
parameter $\d s=e^{\frac{-2\phi(\gamma_r)}{n-m}}\d r$.  
Then 
\begin{align}
\frac{\d\lambda}{\d s}\leq -\frac{\lambda^2}{n-m}-e^{\frac{4\phi(\gamma_r)}{n-m}}\text{
\rm Ric}_{m,n}(L)\left(\dot{\gamma}_r,\dot{\gamma}_r \right).\label{eq:lambda/ds}
\end{align} 
\end{lem}
\begin{pf} We modify the proof of the Bakry-Qian  Laplacian comparison theorem given in Section 10 of \cite{XDLi05}.  Choosing the normal polar coordinate system $(r, \theta)$ at $p\in M$, where $r>0$ and $\theta\in \mathbb{S}^{n-1}$. Let  $J_{\phi}=e^{-\phi}\sqrt{{\rm det}g}$. 
Denote ${'}={\partial \over \partial r}$ and $''={\partial ^2\over \partial r^2}$.  In p.~1355 (see line 5 from the bottom) in \cite{XDLi05}, the following identity has been proved 
\begin{align}
{J_\phi'' \over J_\phi} =-\sum\limits_{i, j} h_{ij}^2-\Ric_{\infty,n}(L) \left({\partial \over \partial r}, {\partial \over \partial r}\right)+(H-\phi')^2, \label{JHR}
\end{align}
where $h_{ij}$ denotes the second fundamental form of $\partial B_r(p)$ at $x=(r, \theta)$ with respect to the unit normal vector ${\partial \over \partial r}$, and $H=\sum\limits_{i}h_{ii}$. 

Let $u={J_\phi'  \over J_\phi} $.  By $(9.53)$ in p. 1353 in \cite{XDLi05}, we have
\begin{eqnarray}
u=Lr=\Delta r-\phi'=H-\phi'.\label{u=Lr}
\end{eqnarray}
Combining $(\ref{JHR})$ and $(\ref{u=Lr})$, 
we have 

\begin{eqnarray}
u'=
 -\sum\limits_{i, j} h_{ij}^2-\Ric_{\infty,n}(L) \left({\partial \over \partial r}, {\partial \over \partial r}\right).\label{JHR2}
\end{eqnarray}
Notice that 
\begin{align*}
\sum\limits_{i, j}h_{ij}^2\geq \sum\limits_{i}h_{ii}^2\geq {\left(\sum\limits_{i=1}^{n-1} h_{ii} \right)^2\over n-1}={H^2\over n-1}={(\Delta r)^2\over n-1}.
\end{align*}
Since $m<1$,  we have the following inequality
\begin{align}
{\d\over \d r} Lr  \leq -{(\Delta r)^2\over n-m}
-\Ric_{\infty,n}(L)
\left( {\partial \over \partial r},  
{\partial \over \partial r}
\right).
\label{Lapl}
\end{align}
This gives us the following inequality  along $\gamma$, 
\begin{align}
\frac{\d}{\d r}(Lr_p)(r,\theta)&\leq -\frac{(\Delta r_p(r,\theta))^2}{n-m}-{\Ric}_{m,n}(L)
\left(\dot{\gamma}_r,\dot{\gamma}_r \right)+\frac{1}{n-m}|\langle  \nabla \phi,\nabla r_p\rangle (r,\theta)|^2.\label{eq:equaitonalonggama}
\end{align}

From this, we can deduce 
\begin{align*}
\frac{\d\lambda}{\d s}&=e^{\frac{2\phi(r,\theta)}{n-m}}\frac{\d\lambda}{\d r}
\leq
-\frac{\lambda^2}{n-m}-e^{\frac{4\phi(r,\theta)}{n-m}}{\Ric}_{m,n}(L)\left(\dot{\gamma}_r,\dot{\gamma}_r\right).
\end{align*}
\end{pf}

\begin{remark} 
{\rm Indeed, a variant of the inequality  \eqref{Lapl} has been obtained in \cite{XDLi05}.  In the first line of p.1356 in \cite{XDLi05}, it was proved that for any $m>n$, it holds
\begin{eqnarray}
{J_\phi''\over J_\phi}\leq -\Ric_{m, n}(L)\left({\partial \over \partial r}, {\partial \over \partial r}\right)+{m-n-1\over m-n}\phi'^2+{n-2\over n-1}H^2-2H\phi'.  \label{mHphi}
\end{eqnarray}
Taking $m\rightarrow \infty$ in $(\ref{mHphi})$, we obtain $(\ref{Lapl})$. In  \cite{WylieYeroshkin},  Wylie and Yeroshkin  
proved  \eqref{Lapl}  by a 
different argument and used it to prove a Laplacian comparison theorem for $m=1$ in terms of conformal changing the metric. Lemma 3.1 and Theorem 2.4 extend Wylie and Yeroshkin's Laplacian comparison theorem to the case $m <1$.
}
\end{remark}

Let $\kappa$ be a continuous function on $[0,+\infty[$ with 
respect to the parameter $s$. 
Assuming the curvature bound ${\Ric}_{m,n}(L)_x(\nabla r_p,\nabla r_p)\geq (n-m)\kappa(s_p(x)) {e^{-\frac{4\phi(x)}{n-m}}}$ for  $s_p(x)<S$ with $x\notin {\rm Cut}(p)\cup\{p\}$, we see 
${\Ric}_{m,n}(L)(\dot{\gamma}_r,\dot{\gamma}_r)\geq (n-m)\kappa(s) {e^{-\frac{4\phi(\gamma_r)}{n-m}}}$ for $s=s(r,\theta)<S$ with $0<r<d(p,{\rm Cut}(p))$. 
From \eqref{eq:lambda/ds} we have the usual Riccati inequality
\begin{align}
-\frac{\d\lambda}{\d s}(s)\geq (n-m)\kappa(s)+\frac{\lambda(s)^2}{n-m}\quad\text{ for }\quad s\in]0,S[\label{eq:RiccattiIneq}
\end{align}
with the caveat that it is in terms of the parameter $s$ instead of $r$. 
This gives us the following comparison estimate. 

\begin{lem}
\label{lem:LaplacianComparisonconformal}
Suppose that $(M,g)$ be an $n$-dimensional complete Riemannian manifold and a function $\phi\in C^2(M)$. Take $R\in]0,+\infty[$ and $x,p\in M$. 
Assume that \eqref{eq:RiciLowBdd} holds for $r_p(x)<R$ 
with $x\notin {\rm Cut}(p)\cup\{p\}$.  
Let $\gamma$, $s$, and $\lambda$ be as in Lemma~\ref{lem:RiccatiDiferIneq}. Then 
\begin{align}
\lambda(r,\theta)\leq m_{\kappa}(s) \label{eq:LocalLapComp}
\end{align}
holds for $r<R$ and $s<\delta_{\kappa}$ and $x=(r,\theta)\notin\text{\rm Cut}(p)\cup\{p\}$. 
Here $s=s(r)=\int_0^r\exp\left(-\frac{2\phi(\gamma_t)}{n-m} \right)\d t$. 
\end{lem}
\begin{pf}
The proof is a mimic of the proof of \cite[Lemma~4.2]{WylieYeroshkin} based on \eqref{eq:RiccatiEqM} and \eqref{eq:RiccattiIneq}. We omit the detail. 
\end{pf}

\begin{cor}\label{cor:Cutlocus}
Let $(M,g)$ be an $n$-dimensional complete Riemannian manifold and $\phi\in C^2(M)$. Take $p\in M$ {\color{black}{and $R\in]0,+\infty[$}}. 
Assume that \eqref{eq:RiciLowBdd} holds for $r_p(x)<R$ 
with $x\notin {\rm Cut}(p)\cup\{p\}$.  
Then $s_p(x)<\delta_{\kappa}$ always holds. 
\end{cor}
\begin{pf}
We may assume $\delta_{\kappa}<\infty$. 
Take $x\in B_R(p)$ with $x\notin {\rm Cut}(p)\cup\{p\}$. 
Let $x=(r,\theta)$ be the polar coordinate expression around $p$ and set 
$s:=s(r)= \int_0^{r}\exp\left(-\frac{\;2 \phi(\gamma_t)\;}{n-m} \right)\d t$ and $S:=s(R)$, where $\gamma$ is a unit speed geodesic with 
$\gamma_0=p$ and $\dot{\gamma}_0=\theta$. 
We see $s_p(x)<S$.  
Assume $S>\delta_{\kappa}$. Then 
there exists $r_0\in]0,R[$ such that $\delta_{\kappa}=\int_0^{r_0}\exp\left(-\frac{2\phi(\gamma_t)}{n-m} \right)\d t$. By \eqref{eq:LocalLapComp},   
$\lambda(r,\theta)\leq (n-m)\cot_{\kappa}(s)$ holds for $s<\delta_{\kappa}$. 
Since $r\uparrow r_0$ is equivalent to $s=s(r)\uparrow\delta_{\kappa}$, we have  
$$
\lambda(r_0,\theta)=\lim_{r\uparrow r_0}\lambda(r,\theta)\leq 
\lim_{r\uparrow r_0}(n-m)\cot_{\kappa}(s(r))=-\infty.
$$ This contradicts the 
well-definedness of $\lambda(r,\theta)=(e^{\frac{2\phi}{n-m}}Lr_p)(r,\theta)$ for $r\in]0,R[$. 
Therefore $S\leq\delta_{\kappa}$ under $\delta_{\kappa}<\infty$ and we obtain the conclusion $s_p(x)<S\leq\delta_{\kappa}$. 
\end{pf}

\begin{apf}{Theorem~\ref{thm:GlobalLapComp}} 
The implication \eqref{eq:RiciLowBdd}$\Longrightarrow$\eqref{eq:GloLapComp} 
for $R<\infty$  
follows from Lemma~\ref{lem:LaplacianComparisonconformal}, because 
$r_p$ is smooth on $M\setminus(\text{\rm Cut}(p)\cup\{p\})$. 
The implication \eqref{eq:RiciLowBdd}$\Longrightarrow$\eqref{eq:GloLapComp} 
for $R=+\infty$ follows from it. 
The latter assertion follows from Corollary~\ref{cor:Cutlocus}.
\end{apf}

\section{Proof of Theorem~\ref{thm:WeightedMyers}}

\begin{apf}{Theorem~\ref{thm:WeightedMyers}}
Suppose that there exist points $p,q\in M$ such that $s(p,q)>\delta_{\kappa}$. 
Since $\text{\rm Cut}(p)$ is closed and measure zero, we may assume $q\notin \text{\rm Cut}(p)$. By Lemma~\ref{lem:LaplacianComparisonconformal}, 
along minimal geodesic from $p$ to $q$, 
$\lambda(r,\theta)\leq (n-m)\cot_{\kappa}(s)$. However, as 
$s\to\delta_{\kappa}$, $\cot_{\kappa}(s)\to-\infty$. This implies 
$\Delta r_p(x)\to-\infty$ as $s(p,x)\to\delta_{\kappa}$. This contradicts 
that $r_p$ is smooth in a neighborhood of $q$. Next we suppose the 
$(\phi,m)$-completeness at $p$. 
Suppose that $\sup_{q\in M}d(p,q)=+\infty$. Then there exists a sequence $\{q_i\}$ in $M$ such that $d(p,q_i)\to+\infty$ as $i\to+\infty$. By
Lemma~\ref{lem:phimcomplete}, there exists a subsequence $\{q_{i_k}\}$ so that 
$s(p,q_{i_k})\to+\infty$ as $k\to+\infty$, which contradicts $\sup_{q\in M}s(p,q)\leq\delta_{\kappa}$. 
Therefore, $\sup_{q\in M}d(p,q)<\infty$, hence $M$ is compact. 
\end{apf}

\section{Proof of Theorem~\ref{thm:BGVol} 
 }

\begin{lem}[Volume Element Comparison]\label{lem:VolComp}
Let $(M,g)$ be an $n$-dimensional complete Riemannian manifold and $\phi\in C^2(M)$. Take $p\in M$ {\color{black}{and $R\in]0,+\infty]$}}.  
Assume that \eqref{eq:RiciLowBdd} holds {\color{black}{for $r_p(x)<R$ with  
$x\notin {\rm Cut}(p)\cup\{p\}$}}. 
Let $J$ be the volume element in geodesic polar coordinates {\color{black}{around $p\in M$}} and set $J_{\phi}(r,\theta):=e^{-\phi(r,\theta)}J(r,\theta)$. Then for $r_0<r_1{\color{black}{<R\land {\rm cut}(\theta)}}$,  
\begin{align}
\frac{J_{\phi}(r_1,\theta)}{J_{\phi}(r_0,\theta)}
\leq \frac{\s_{\kappa}(s(r_1,\theta))^{n-m}}{\s_{\kappa}(s(r_0,\theta))^{n-m}}.\label{eq:VolComp} 
\end{align}
\end{lem}
{\color{black}{Here $\text{\rm cut}(\theta)$ is the distance from $p$ to the cut point along the geodesic with $\gamma(0)=p$ and $\dot{\gamma}(0)=\theta$.}}
\begin{pf}
The proof is a mimic of the proof of \cite[Lemma~4.3]{WylieYeroshkin}. 
We omit the detail. 
\end{pf}

\begin{apf}{Theorem~\ref{thm:BGVol}}
Consider geodesic polar coordinates {\color{black}{$(r,\theta)$ around $p$. 
Then }}
\begin{align*}
\mu(A(p,r_0,r_1))=\int_{\mathbb{S}^{n-1}}\int_{{\rm cut}(\theta)\land r_0}^{{\rm cut}(\theta)\land r_1}J_{\phi}(r,\theta)\d r \d \theta.
\end{align*}
By Lemma~\ref{lem:VolComp}, for all $r_1,r_2>0$ with $r_1< r_2<R$ 
{\color{black}{and 
$r_2<{\rm cut}(\theta)$}}
\begin{align*}
\frac{J_{\phi}(r_2,\theta)}{J_{\phi}(r_1,\theta)}\leq \frac{\s_{\kappa}^{n-m}(s_p(r_2,\theta))}{\s_{\kappa}^{n-m}(s_p(r_1,\theta))}\leq 
\frac{\s_{\kappa}^{n-m}\left(
{\color{black}{\sup_{\eta\in\mathbb{S}^{n-1}}s_p(r_2,\eta)}}
\right)}{\s_{\kappa}^{n-m}\left(
{\color{black}{\inf_{\eta\in\mathbb{S}^{n-1}}s_p(r_1,\eta)}}
\right)}.
\end{align*}
So for 
{\color{black}{$0\leq r_a< r_b\leq r_d$, $0\leq r_a\leq r_c< r_d$ and $r_d<R$}}, 
we have {\color{black}{the following inequality}}
\begin{align*}
\frac{\int_{\text{\rm cut}(\theta)\land r_c}^{\text{\rm cut}(\theta)\land r_d}J_{\phi}(r_2,\theta)\d r_2}{\int_{\text{\rm cut}(\theta)\land r_a}^{\text{\rm cut}(\theta)\land r_b}J_{\phi}(r_1,\theta)\d r_1}
&
\leq\frac{\int_{\text{\rm cut}(\theta)\land r_c}^{\text{\rm cut}(\theta)\land r_d}\s_{\kappa}^{n-m}(s_p(r_2,\theta))\d r_2}{\int_{\text{\rm cut}(\theta)\land r_a}^{\text{\rm cut}(\theta)\land r_b}\s_{\kappa}^{n-m}(s_p(r_1,\theta))\d r_1}
\\
&\leq 
\frac{\int_{r_c}^{r_d}\s_{\kappa}^{n-m}\left(
{\color{black}{\sup_{\eta\in\mathbb{S}^{n-1}}s_p(r_2,\eta)}}
\right)\d r_2}{\int_{ r_a}^{ r_b}\s_{\kappa}^{n-m}\left(
{\color{black}{\inf_{\eta\in\mathbb{S}^{n-1}}s_p(r_1,\eta)}}
\right)\d r_1}
\end{align*}
{\color{black}{under $r_a=r_c$ or $r_b=r_d$ by use of \cite[Lemma~3.1]{Zhu97} 
(cf.~\cite[Proof of Theorem~3.2]{Zhu97})}}.
Thus {\color{black}{we can deduce}} 
\begin{align*}
\frac{\int_{\mathbb{S}^{n-1}}\int_{\text{\rm cut}(\theta)\land r_c}^{\text{\rm cut}(\theta)\land r_d}J_{\phi}(r_2,\theta)\d r_2\d\theta}{\int_{\mathbb{S}^{n-1}}\int_{\text{\rm cut}(\theta)\land r_a}^{\text{\rm cut}(\theta)\land r_b}J_{\phi}(r_1,\theta)\d r_1\d\theta}
\leq\frac{\int_{\mathbb{S}^{n-1}}\int_{r_c}^{r_d}\s_{\kappa}^{n-m}\left(
{\color{black}{\sup_{\eta\in\mathbb{S}^{n-1}}s_p(r_2,\eta)}}
\right)\d r_2\d\theta}{\int_{\mathbb{S}^{n-1}}\int_{r_a}^{r_b}\s_{\kappa}^{n-m}\left(
{\color{black}{\inf_{\eta\in\mathbb{S}^{n-1}}s_p(r_1,\eta)}}
\right)\d r_1\d\theta}
\end{align*}
{\color{black}{for general $0\leq r_a< r_b\leq r_d$ and $0\leq r_a\leq r_c< r_d$}}.
This implies {\color{black}{that}} \eqref{eq:BGAnnuliUpLow} {\color{black}{holds for $r_1<R$}}. 
\end{apf}

\section{Proof of Theorem~\ref{thm:conservative}}  
The proof of Theorem~\ref{thm:conservative} is based on the generalized Grigoryan's criterion for the conservativeness of Dirichlet form, which says that if for some $x\in M$, 
\begin{align}
\int_1^{\infty}\frac{r\,\d r}{\log \mu(B_r(p))}=+\infty,\label{eq:Grigoryan}
\end{align}
then $P_t=e^{tL}$ is conservative. In view of this, Theorem~\ref{thm:conservative} follows from the following lemma.
\begin{lem}\label{lem:VolumeComparison}
Take $p\in M$ {\color{black}{and $R\in]0,+\infty]$}}. 
Let $(M,g)$ be an $n$-dimensional complete Riemannian manifold without boundary and a function $\phi\in C^2(M)$. 
Let ${\sf K}(s)$ be a non-negative continuous non-decreasing  function on $[0,+\infty[$ satisfying 
${\bf K}(p)$.  
 Assume that \eqref{eq:RiciLowBdd} holds {\color{black}{for $r_p(x)<R$ with }} 
$x\notin {\rm Cut}(p)\cup\{p\}$.
Then, {\color{black}{for 
$0<r_1<r_2<R$}}
\begin{align*}
{\color{black}{\frac{\mu(B_{r_2}(p))}{\mu(B_{r_1}(p))}}}
&
{\color{black}{ \leq  
e^{2(\overline{\phi}_p(r_1)-\underline{\phi}_p(r_2))}
\left(\frac{r_2}{r_1}\right)^
{\left(n-m+1\right)}
\exp\left(\sqrt{(n-m){\sf K}\left(e^{-\frac{2\underline{\phi}_p(r_2)}{n-m}}r_2 \right)}
e^{-\frac{2\underline{\phi}_p(r_2)}{n-m}}r_2 \right).}}
\end{align*}
In particular, under  $\text{\rm Ric}_{m,n}(L)_x(\nabla r_p,\nabla r_p)\geq0$ for 
$x\notin{\rm Cut}(p)\cup\{p\}$, we have that {\color{black}{for $0<r_1<r_2<R$}} 
\begin{align}
{\color{black}{
\frac{\mu(B_{r_2}(p))}{\mu(B_{r_1}(p))}\leq
 e^{2(\overline{\phi}_p(r_1)-\underline{\phi}_p(r_2))}
\left(\frac{r_2}{r_1}\right)^{\left(n-m+1\right)}}}.\label{eq:BGBall}
\end{align}
\end{lem}
\begin{pf}
Recall that ${\sf K}(s)=(n-m)\kappa(s)$ is a non-negative non-decreasing continuous function. 
Assume that 
$\text{\rm Ric}_{m,n}(L)_x(\nabla r_p,\nabla r_p)\geq -K(x)e^{-\frac{4\phi(x)}{n-m}}$ holds for {\color{black}{$r_p(x)<R$ with}} $x\in 
 ({\rm Cut}(p)\cup\{p\})^c$. Here $K(x)={\sf K}(s_p(x))$, $x\in M$. 
By Theorem~\ref{thm:BGVol}, 
we have {\color{black}{for $0<r_1<r_2<R$
\begin{align}
\frac{\mu(B_{r_2}(p))}{\mu(B_{r_1}(p))}\leq
\frac{\int_0^{r_2}s_{-\kappa}^{n-m}\left(
e^{-\frac{2\underline{\phi}_p(r)}{n-m}}
r\right)\d r}
{\int_0^{r_1} s_{-\kappa}^{n-m}\left(e^{-\frac{2\overline{\phi}_p(r)}{n-m}}r\right)\d r},\label{eq:BishopGromov}
\end{align} 
where 
$s_{-\kappa}(s)$ is the unique solution of the Jacobi equation $s_{-\kappa}''(s)-\kappa(s)s_{-\kappa}(s)=0$ with $s_{-\kappa}(0)=0$ and $s_{-\kappa}'(0)=1$. 
Let $T:=e^{-\frac{2\underline{\phi}_p(r_2)}{n-m}}r_2$. 
Applying the comparison theorem for the solution of 
Jacobi equation, we see 
\begin{align}
s=s_0(s)\leq s_{-\kappa}(s)\leq s_{-\kappa(T)}(s)=\frac{\sinh \sqrt{\kappa(T)}s}{\sqrt{\kappa(T)}}. \label{eq:JacobiComp}
\end{align}
Combining \eqref{eq:BishopGromov} 
with \eqref{eq:JacobiComp}, 
we have 
\begin{align*}
\frac{\mu(B_{r_2}(p))}{\mu(B_{r_1}(p))}&\leq \frac{e^{2\overline{\phi}_p(r_1)}}{\int_0^{r_1} r^{n-m}\d r}\int_0^{r_2} s_{-\kappa(T)}^{n-m}\left(e^{-\frac{2\underline{\phi}_p(r)}{n-m}}r \right)\d r\\
&\leq 
\frac{e^{2\overline{\phi}_p(r_1)}}{\int_0^{r_1} r^{n-m}\d r}\int_0^{r_2}
\left(\frac{\sinh \sqrt{\kappa(T)}
e^{-\frac{2\underline{\phi}_p(r)}{n-m}}r
}{\sqrt{\kappa(T)}}\right)^{n-m}
\d r\\
&\leq 
\frac{e^{2\overline{\phi}_p(r_1)}}{\int_0^{r_1} r^{n-m}\d r}\int_0^{r_2}
 \left(\frac{\sinh \sqrt{\kappa(T)}T
}{\sqrt{\kappa(T)}T}\cdot 
e^{-\frac{2\underline{\phi}_p(r_2)}{n-m}}
 r\right)^{n-m}
\d r,
\end{align*}
}}
where we use that $x\mapsto \sinh x/x$ is non-decreasing. 
Thus
{\color{black}{
\begin{align*}
\frac{\mu(B_{r_2}(p))}{\mu(B_{r_1}(p))}\leq 
e^{2(\overline{\phi}_p(r_1)-\underline{\phi}_p(r_2))}\left(\frac{r_2}{r_1} \right)^{n-m+1}
\left(\frac{\sinh \sqrt{\kappa(T)}T}{\sqrt{\kappa(T)}T} \right)^{n-m}.
\end{align*}
}}
Using the inequality $\frac{\sinh x}{x}\leq e^x$ for $x\geq0$, we obtain the conclusion. \end{pf}
\begin{apf}{Theorem~\ref{thm:conservative}}
We first prove the assertion 
 for the case ${\sf K}\equiv0$ under 
\eqref{eq:asymptoticphi}. {\color{black}{Let $r_0$ be the constant specified in 
\eqref{eq:asymptoticphi}. }}  
From \eqref{eq:BGBall}, we have that for $r>r_0{\color{black}{>0}}$
\begin{align*}
\log\mu(B_r(p))&\leq 
\log\mu(B_{r_0}(p))+2\overline{\phi}_p(r_0)-2\underline{\phi}_p(r)+
(n-m+1)\log{\color{black}{\left(r/r_0\right)}}.
\end{align*}
There exists $r_1>r_0$ such that for all $r>r_1$
\begin{align*}
\log \mu(B_{r_0}(p))+2\overline{\phi}_p(r_0){\color{black}{-(n-m+1)\log r_0}}\leq (n-m+1)\log r. 
\end{align*}
Thus, for all $r>r_1$
\begin{align*}
\log\mu(B_r(p))\leq 2(n-m+1)\log r-2\underline{\phi}_p(r).
\end{align*}
Since $\lim_{r\to\infty}\frac{\log r}{r}=0$ and $\lim_{r\to\infty}\exp\left(-\frac{2\underline{\phi}_p(r)}{n-m} \right)=C\in ]0,+\infty]$, there exists $r_2\in]r_1,+\infty[$ such that for any $r>r_2$  
\begin{align*}
2(n-m+1)\frac{\log r}{r}\leq \frac{n-m}{r_0}\exp\left({-\frac{2\underline{\phi}_p(r)}{n-m}}\right).
\end{align*}
From this 
\begin{align*}
\frac{\log\mu(B_r(p))}{r}&\leq  
2(n-m+1)\frac{\log r}{r}+\frac{n-m}{r}\left(-\frac{2\underline{\phi}_p(r)}{n-m}\right)
\leq 2\frac{n-m}{r_0}\exp\left(-\frac{2\underline{\phi}_p(r)}{n-m}\right).
\end{align*}
By \eqref{eq:asymptoticphi}, we have  
\begin{align*}
\int_{r_2}^{\infty}\frac{r\d r}{\log\mu(B_r(p))}\geq \frac{r_0}{2(n-m)}
\int_{r_2}^{\infty}\exp\left(\frac{2\underline{\phi}_p(r)}{n-m}\right)\d r
=+\infty,
\end{align*}
which implies the conservativeness of {\bf X} by \cite{Grigo:Conservative}. 
When $n\leq m+1$ and $\phi$ is lower bounded, we see 
\begin{align*}
\int_1^{\infty}\frac{r\,\d r}{\mu(B_r(p))}\geq 
\frac{e^{2\inf_M\phi-2\overline{\phi}_p(r_0)}r_0^{n-m+1}}{\mu(B_{r_0}(p))}\int_{r_1}^{\infty}\frac{\d \,r}{r^{n-m}}=+\infty.
\end{align*}
This implies the recurrence of {\bf X} (see \cite[Theorem~3]{St:DirI}).
Next we prove the assertion for the case ${\sf K}\not\equiv0$. 
Then there exists $t_0>0$ such that {\color{black}{
${\sf K}(t)\geq {\sf K}(t_0)>0$}} for all $t\geq t_0$. 
{\color{black}{
Let $r_0$ be the constant specified in \eqref{eq:integrability}. 
Replacing ${\sf K}$ with ${\sf K}\lor {\sf K}(t_0)$, 
we may assume that ${\sf K}$ is bounded below by a positive constant.}}
Note that $r\mapsto e^{-\frac{2\underline{\phi}_p(r)}{n-m}}$ is lower bounded by $e^{-\frac{2{\phi}(p)}{n-m}}$.  
In this case, there exists {\color{black}{$r_1>r_0$ such that   
\begin{align*}
\frac12\sqrt{(n-m){\sf K}\left(e^{-\frac{2\underline{\phi}_p(r)}{n-m}}r\right)}e^{-\frac{2\underline{\phi}_p(r)}{n-m}}r
\geq\log \mu(B_{r_0}(p))
 +2\overline{\phi}_p(r_0)
+(n-m+1)\log r\quad\text{ for all }\quad r>r_1.
\end{align*}
}} 
If we set 
\begin{align*}
{\color{black}{r_2}}:=\max\left\{\frac{2}{\sqrt{(n-m){\sf K}\left(e^{\frac{-2\underline{\phi}_p(r_0)}{n-m}}r_0 \right)}},{\color{black}{r_1}}\right\},
\end{align*} 
then 
\begin{align*}
\frac12\sqrt{(n-m){\sf K}\left(e^{-\frac{2\underline{\phi}_p(r)}{n-m}}r\right)}e^{-\frac{2\underline{\phi}_p(r)}{n-m}}r\geq -2\underline{\phi}_p(r) 
\quad\text{ for all }\quad r>{\color{black}{r_2}}.
\end{align*}
Thus
\begin{align*}
\sqrt{(n-m){\sf K}\left(e^{-\frac{2\underline{\phi}_p(r)}{n-m}}r\right)}e^{-\frac{2\underline{\phi}_p(r)}{n-m}}r\geq\log \mu(B_{r_0}(p))
 +2\overline{\phi}_p(r_0)-2\underline{\phi}_p(r)
+(n-m+1)\log r
\end{align*}
for all $r>{\color{black}{r_2}}$.   
Applying this with Lemma~\ref{lem:VolumeComparison}, we have 
\begin{align*}
\int_{{\color{black}{r_2}}}^{\infty}&\frac{r\d r}{\log \mu(B_r(p))}\\&
\hspace{-0.4cm}
\geq \int_{{\color{black}{r_2}}}^{\infty}\!\!\frac{r\d r}{\log \mu(B_{r_0}(p))
 \!+\!2\overline{\phi}_p(r_0)\!-\!2\underline{\phi}_p(r)\!+\!
(n-m+1)\log r\!+\!\sqrt{(n\!-\!m){\sf K}\left(\!e^{-\frac{2\underline{\phi}_p(r)}{n-m}}r\!\right)}e^{-\frac{2\underline{\phi}_p(r)}{n-m}}r}\\
&\hspace{-0.4cm}\geq \frac12\int_{{\color{black}{r_2}}}^{\infty}\frac{r\d r}{\sqrt{(n-m){\sf K}\left(\!e^{-\frac{2\underline{\phi}_p(r)}{n-m}}r\!\right)}e^{-\frac{2\underline{\phi}_p(r)}{n-m}}r}
=\frac{1}{2\sqrt{n-m}}\int_{r_1}^{\infty}\frac{\d r}{\sqrt{{\sf K}\left(e^{-\frac{2\underline{\phi}_p(r)}{n-m}}r\right)}e^{-\frac{2\underline{\phi}_p(r)}{n-m}}}=+\infty.
\end{align*}
\end{apf}

\section{Proof of Theorem~\ref{thm:Fellerproperty}}

\begin{apf}{Theorem~\ref{thm:Fellerproperty}}
 We follow the argument as used in the proof of Theorem 1.5 in \cite{XDLi05}, which extends the method originally by Azencott~\cite{Azencott:Behavi}    
and developed in Hsu~\cite[Theorem~4.3.2]{Hsu:2001} (see also Hsu~\cite{Hsu:1989} and 
Qian~\cite{Qian:ConservativeFeller}). 
By Theorem~\ref{thm:conservative}, the $L$-diffusion is 
conservative under the condition of Theorem~\ref{thm:Fellerproperty}. 
Let ${\bf X}=(\Omega, X_t, \P_x)$ be the $L$-diffusion starting from $x\in M$. 
By Azencott~\cite{Azencott:Behavi}, we need to prove that for any geodesic ball 
$K=B_R(p)$, where $R$ is a fixed constant, we have 
\begin{align}
\lim_{d(x,p)\to\infty}\P_x(\sigma_K<t)=0,\label{eq:CriterionFeller}
\end{align}    
where $\sigma_K:=\inf\{t>0\mid X_t\in K\}$ is the first hitting time to $K=B_R(p)$. 
We may assume $R<r_p(x)$. 
Let $\sigma_0:=0$, and for all $k\in\mathbb{N}$, 
\begin{align*}
\tau_k:=&\inf\{t>\sigma_k\mid d(X_t,X_{\sigma_k})=1\},\quad k\geq0,\\
\sigma_k:=&\inf\{t\geq\tau_{k-1}\mid r_p(X_t)=r_p(x)-k\},\quad k\geq1.
\end{align*}
That is, $\sigma_k$ is the first hitting time to the geodesic ball 
$B_{r_p(x)-k}(p)$, $\tau_k-\sigma_k$ is the amount of time during which the $L$-diffusion process moves from $X_{\sigma_k}\in \partial B_{r_p(x)-k}(p)$ to 
$X_{\tau_k}\in \partial B_1(X_{\sigma_k})$, and $\sigma_{k+1}-\tau_k$ is the amount of the time during which the $L$-diffusion leaves from $\partial B_1(X_{\sigma_k})$ and hits $\partial B_{r_p(x)-k+1}(p)$. Let  
\begin{align*}
\theta_k:=\tau_k-\sigma_k.
\end{align*}
Then 
\begin{align*}
\sigma_K\geq\sigma_{\lfloor r_p(x)-R\rfloor }\geq\theta_0+\theta_1+\cdots+\theta_{\lfloor r_p(x)-R-1\rfloor },
\end{align*}
where $\lfloor r_p(x)-R\rfloor$ denotes the largest integer which does not exceed $r_p(x)-R$.  
Since ${\sf K}(s)=(n-m)\kappa(s)$ is non-decreasing, we may assume 
that there exists $t_0>0$ satisfying $\kappa(t)>0$ for all 
$t\geq t_0$ provided ${\sf K}\not\equiv0$. 
The key point is to prove that there exist two constants $C_1,C_2>0$ such that for all $k\geq0$, 
\begin{align}
\P_x&\left(\theta_k\leq \frac{C_1}{
l(r_p(x)-k+1)} \right)\leq \exp\left(-C_2
l(r_p(x)-k+1)\right),\label{eq;iterationest}
\end{align}
where  $l(s):=e^{\frac{2\overline{\phi}_p(s)-2\underline{\phi}_p(s)}{n-m}}\sqrt{\kappa\left(e^{-\frac{2\underline{\phi}_p(s)}{n-m}}s\right)}$ under 
${\sf K}\not\equiv0$, or 
$l(s):=e^{\frac{2\overline{\phi}_p(s)-2\underline{\phi}_p(s)}{n-m}}$ under ${\sf K}\equiv0$.  
To this end, we use Kendall's It\^o-Skorokhod formula. In fact, see Kendall~\cite{Kendall}, under the probability measure $\P_x$, there exists a standard Brownian motion 
$\beta_t$ such that $r_x(X_t)=d(X_t,x)$ can be decomposed into 
\begin{align*}
r_x(X_t)=\sqrt{2}\beta_t+\int_0^t Lr_x(X_s)\1_{\{X_s\notin{\rm Cut}(x)\}}\d s-L_t,\qquad t\in[0,+\infty[,
\end{align*}
where $L_t$ is a non-decreasing process which is increasing only on 
$\{t\in[0,+\infty[\;\mid X_t\in \text{\rm Cut}(x)\}$. 
For a proof, see \cite[Remark~4.1]{XDLi05}. 
Note here that {\bf X} is conservative under \eqref{eq:RicciBound} 
by Theorem~\ref{thm:conservative} with Remark~\ref{rem:Integrability}. 
Moreover, using the Kendall's decomposition and the Girsanov transform, we have (cf.~Qian~\cite{Qian:ConservativeFeller} and \cite{XDLi05} p.1320)
\begin{align*}
d(X_t,X_{\sigma_k})=\sqrt{2}(\beta_t-\beta_{\sigma_k})+\int_{\sigma_k}^tLd(X_s,X_{\sigma_k})\1_{\{X_s\notin{\rm Cut}(X_{\sigma_k})\}}\d s-(L_t-L_{\sigma_k}).
\end{align*}
Note that 
\begin{align*}
d^2(X_t,X_{\sigma_k})=2\int_{\sigma_k}^td(X_s,X_{\sigma_k})\d (d(X_s,X_{\sigma_k}))
+\langle  d(X_{\cdot},X_{\sigma_k})\rangle _t-\langle  d(X_{\cdot},X_{\sigma_k})\rangle _{\sigma_k}.
\end{align*}
Since $\langle  d(X_{\cdot},X_{\sigma_k})\rangle _t=\langle  \sqrt{2}\beta\rangle _t=2t$ and $L_t-L_{\sigma_k}$ 
is a non-decreasing positive process on $[\sigma_k,\tau_k]$, we have 
\begin{align}
\frac12d^2(X_t,X_{\sigma_k})\leq \sqrt{2}\int_{\sigma_k}^td(X_s,X_{\sigma_k})\d \beta_s
+\int_{\sigma_k}^td(X_s,X_{\sigma_k})Ld(X_s,X_{\sigma_k})
\1_{\{X_s\notin{\rm Cut}(X_{\sigma_k})\}}
\d s+t-\sigma_k.\label{eq:comparisonsquare}
\end{align}
For $s\in[\sigma_k,\tau_k]$, $X_s\in \overline{B_1}(X_{\sigma_k})\subset B_{r_p(x)-k+1}(p)$. 
This means that $y:=X_s$ and $q:=X_{\sigma_k}$ satisfy $r_p(y)<r_p(x)-k+1$. 
Applying \eqref{eq:spIneQrp} to $y$, we can see that 
\begin{align}
e^{-\frac{2\overline{\phi}_q(r_q(y))}{n-m}}r_q(y)\leq s_q(y)
\quad \text{ and }\quad 
s_p(y)
<
e^{-\frac{2\underline{\phi}_p(r_p(x)-k+1)}{n-m}}(r_p(x)-k+1).\label{eq:sp}
\end{align} 
From \eqref{eq:RicciBound}, we have 
\begin{align}
{\Ric}_{m,n}(L)(y)\geq -{\sf K}(s_p(y))e^{-\frac{4\phi(y)}{n-m}}\label{eq:Riccisp}
\end{align}
under $r_q(y)\leq1$. 
Now we apply Theorem~\ref{thm:GlobalLapComp} to \eqref{eq:Riccisp}
 as  $q=X_{\sigma_k}$ is the reference point and $s=s_p(y)$ is a constant. Then, 
for $y\in \overline{B_1}(q)\setminus(\text{\rm Cut}(q)\cup\{q\})$ 
\begin{align}
(Lr_q)(y)&\leq 
(n-m)\sqrt{\kappa(s_p(y))}\coth\left(
\sqrt{\kappa(s_p(y))}
s_q(y) \right)e^{-\frac{2\phi(y)}{n-m}}\notag\\
&
\leq 
(n-m)\sqrt{\kappa(s_p(y))}\coth\left(
\sqrt{\kappa(s_p(y))}
e^{-\frac{2\overline{\phi}_q(r_q(y))}{n-m}} r_q(y)\right)e^{-\frac{2\phi(y)}{n-m}}
\notag\\
&\leq 
\frac{n-m}{r_q(y)}e^{\frac{2\overline{\phi}_q(r_q(y))}{n-m}}\left(1+\sqrt{\kappa(s_p(y))}
e^{-\frac{2\overline{\phi}_q(r_q(y))}{n-m}} r_q(y) \right)e^{-\frac{2\phi(y)}{n-m}},
\label{eq:Est}
\end{align}
where we use the elementary inequality $a\coth a\leq 1+a$ for $a\geq0$. 
Applying \eqref{eq:spIneQrp}, $r_q(y)\leq1$ and $r_p(q)=r_p(x)-k$ together imply $\overline{\phi}_q(r_q(y))\leq \overline{\phi}_p(r_q(y)+r_p(q))\leq \overline{\phi}_p(r_p(x)-k+1)$ 
and $\phi(y)\geq\underline{\phi}_q(r_q(y))\geq \underline{\phi}_p(r_q(y)+r_p(q))\geq \underline{\phi}_p(r_p(x)-k+1)$. 
Combining these with 
\eqref{eq:sp}, 
for $q=X_{\sigma_k}$ and $y=X_s\notin {\rm Cut}(X_{\sigma_k})$ under $s\in[\sigma_k,\tau_k]$, 
we obtain 
\begin{align*}
d(X_s,X_{\sigma_k})Ld(X_s,X_{\sigma_k})&\leq (n-m)e^{\frac{2(\overline{\phi}_p-\underline{\phi}_p)(r_p(x)-k+1)}{n-m}}\nonumber\\
&\hspace{-1.4cm}
+(n-m)e^{\frac{2(\overline{\phi}_p-\underline{\phi}_p)(r_p(x)-k+1)}{n-m}}
e^{-\frac{2\phi(p)}{n-m}}
\sqrt{\kappa\left(e^{-\frac{2\underline{\phi}_p(r_p(x)-k+1)}{n-m}}(r_p(x)-k+1) \right)}.
\end{align*} 
When ${\sf K}\not\equiv0$, 
without loss of the generality we may assume 
 \begin{align*}
 \kappa\left(e^{-\frac{2\underline{\phi}_p(r_p(x)-k+1)}{n-m}}(r_p(x)-k+1) \right)
 \geq 1
\end{align*} 
by changing $\kappa(t)$ into $
\kappa(t)/\kappa(t_0)$ and taking sufficiently large $r_p(x)$. Hence, whether ${\sf K}\not\equiv0$ or not, we may assume 
\begin{align*}
l(r_p(x)-k+1)\geq 1.
\end{align*}
Taking $t=\tau_k$ in \eqref{eq:comparisonsquare} and since $d(X_s,X_{\sigma_k})\leq 
d(X_{\tau_k},X_{\sigma_k})=1$ for all $s\in[\sigma_k,\tau_k]$, we obtain 
\begin{align*}
\frac12\leq \sqrt{2}\int_{\sigma_k}^{\tau_k}d(X_s,X_{\sigma_k})\d \beta_s+
2\max\{1,e^{-\frac{2\phi(p)}{n-m}}\}(n-m+1)l(r_p(x)-k+1)(\tau_k-\sigma_k).
\end{align*}
This yields that, for any enough small constant $C_1>0$,
\begin{align*}
\P_x\left(\tau_k-\sigma_k\leq C_1\Bigl/l(r_p(x)-k+1) \right)
\leq\P_x\left(\int_{\sigma_k}^{\tau_k}d(X_s,X_{\sigma_k})\d\beta_s\geq\frac18\right).
\end{align*}
Based on L\'evy's criterion and the random time change, the standard method as used in Hsu~\cite[Theorem~3.6.1]{Hsu:2001}, \cite[Lemma~3.2]{Hsu:1989} proves that 
\begin{align*}
\P_x\left(\int_{\sigma_k}^{\tau_k}d(X_s,X_{\sigma_k})\d\beta_s\geq\frac18\right)
\leq \exp\left(-C_2l(r_p(x)-k+1) \right).
\end{align*}
Therefore we have proved \eqref{eq;iterationest}. Then we can follow the same argument used in Hsu~\cite{Hsu:2001}, \cite{Hsu:1989} to obtain 
\begin{align}
\P_x(\sigma_K\leq t)\leq \sum_{k=0}^{N(x,t)}\exp\left(-C_2l(r_p(x)-k+1) \right),
\label{eq:SuffFeller}
\end{align}
where $N(x,t)$ is the smallest integer such that 
\begin{align*}
\sum_{k=0}^{N(x,t)}\frac{1}{l(r_p(x)-k+1)}>\frac{t}{C_1}.
\end{align*}
Indeed, if ${\sf K}\not\equiv0$ (resp.~${\sf K}\equiv0$), by \eqref{eq:upintegrability} (resp.~\eqref{eq:integarbilityK0}), we can deduce 
\begin{align*}
\int_{s_o}^{\infty}\frac{\d s}{l(s)}
=+\infty
\end{align*}
for some $s_o>0$.  
Such $N(x,t)$ exists for all sufficiently large $r_p(x)$. 
By the choice of 
$N(x,t)$, 
\begin{align}
\frac{t}{C_1}\geq \sum_{k=0}^{N(x,t)-1}\frac{1}{l(r_p(x)-k+1)}\geq \sum_{j=\lfloor r_p(x)\rfloor-N(x,t)+3}^{\lfloor r_p(x)\rfloor +2}\frac{1}{l(j)}\geq \int_{\lfloor r_p(x)\rfloor -N(x,t)+3}^{\lfloor r_p(x)\rfloor +3}\frac{\d r}{l(r)}\label{eq:tC1}
\end{align}
By \eqref{eq:integrability} again, as $r_p(x)\uparrow+\infty$, the lowest bound $\lfloor r_p(x)\rfloor -N(x,t)+3$ of the last sum must go to infinity:
\begin{align}
r_p(x)-N(x,t)\to\infty\quad\text{ as }\quad r_p(x)\uparrow+\infty.\label{eq:divergence}
\end{align}
This implies that $\lfloor r_p(x)-R-1\rfloor \geq N(x,t)$ for all sufficiently large $r_p(x)$, and the following sequence of inclusions holds:
\begin{align*}
\{\sigma_K\leq t\}&\subset \left\{\sum_{k=0}^{\lfloor r_p(x)-R-1\rfloor }\theta_k\leq t \right\}\subset \left\{\sum_{k=0}^{N(x,t)}\theta_k\leq t \right\}\subset \bigcup_{k=0}^{N(x,t)}\left\{\theta_k\leq \frac{C_1}{l(r_p(x)-k+1)} \right\}.
\end{align*}
 Then we can get \eqref{eq:SuffFeller} by \eqref{eq;iterationest}. 
Combining \eqref{eq:tC1} with 
\begin{align*}
e^{-C_2l(r)}\leq l(r_p(x)-N(x,t))e^{-C_2l(r_p(x)-N(x,t))}\frac{1}{l(r)},
\end{align*}
for $r\geq r_p(x)-N(x,t)$,
we obtain 
\begin{align*}
\P_x(\sigma_K\leq t)&\leq \int_{r_p(x)-N(x,t)}^{r_p(x)+1}e^{-C_2l(r)}\d r\\
&\leq l(r_p(x)-N(x,t))e^{-C_2 l(r_p(x)-N(x,t))}\left(\frac{t}{C_1}+\frac{3}{l(r_p(x)-N(x,t))} \right).
\end{align*} 
We may assume $\lim_{s\to\infty}l(s)=+\infty$. 
Indeed, if $\phi$ is unbounded, then this holds automatically. If $\phi$ is bounded (in this case 
the conditions ${\bf K}(p)$ and $\overline{\bf K}(p)$ are equivalent to \eqref{eq:Inttegrability}), 
this holds under 
$\lim_{r\to\infty}{\sf K}(r)=+\infty$. When 
$\lim_{r\to\infty}{\sf K}(r)<+\infty$ including the case 
${\sf K}\equiv0$, we can replace 
${\sf K}(r)$ with another continuous strictly increasing function 
$\widetilde{\sf K}(r):={\sf K}(r)+r$. Then we see
\begin{align*}
{\sf K}(r)\leq\widetilde{\sf K}(r)\quad \text{ for all }r>0,\quad \lim_{r\to\infty}\widetilde{\sf K}(r)=+\infty,\quad\text{ and }\quad\widetilde{\sf K}\quad \text{ satisfies }\quad\eqref{eq:Inttegrability}. 
\end{align*}
The last term  above converges to $0$ 
as $ r_p(x)\to\infty$ by $\lim_{s\to\infty}l(s)=+\infty$ and 
\eqref{eq:divergence}. 
This proves the desired key estimate \eqref{eq:CriterionFeller} for $K={B_R}(p)$. 
\end{apf}
\noindent{\bf Acknowledgements.} 
The authors would like to thank the anonymous referee. His/Her 
comments help us to improve the content of this paper.  

\providecommand{\bysame}{\leavevmode\hbox to3em{\hrulefill}\thinspace}
\providecommand{\MR}{\relax\ifhmode\unskip\space\fi MR }
\providecommand{\MRhref}[2]{%
  \href{http://www.ams.org/mathscinet-getitem?mr=#1}{#2}
}
\providecommand{\href}[2]{#2}

\end{document}